\documentclass[a4paper]{article}
\usepackage{geometry}
\usepackage{graphicx}
\usepackage{amsmath}
\usepackage{amsfonts}
\usepackage{amssymb}

\usepackage{dsfont,xfrac,sidecap,caption}

\usepackage[ruled, vlined]{algorithm2e}

\usepackage{cleveref}
\usepackage{todonotes}

\usepackage{multirow}

\usepackage{pgfplots}

\usepackage{todonotes}

\newcommand{\nt}{\vec{n}}

\newcommand{\vt}{\mathbf{v}}

\newcommand{\ut}{\mathbf{u}}
\newcommand{\wt}{\mathbf{w}}
\newcommand{\sigmat}{\boldsymbol{\sigma}}
\newcommand{\Sigmat}{\boldsymbol{\Sigma}}

\newcommand{\Ft}{\mathbf{F}}

\newcommand{\bs}[1]{\boldsymbol{#1}}

\newcommand{\FL}{{\cal F}}
\newcommand{\SO}{{\cal S}}
\newcommand{\IN}{\Gamma}

\begin{document}

\title{Efficient coarse correction for parallel time-stepping in plaque growth simulations}

\date{}

\author{Stefan Frei\thanks{Department of Mathematics and Statistics, University of Konstanz, 78457 Konstanz, Germany, stefan.frei@uni-konstanz.de},\, Alexander Heinlein\thanks{Delft University of Technology, Delft Institute of Applied Mathematics, Mekelweg 4, 2628CD Delft, Netherlands, a.heinlein@tudelft.nl}}

\maketitle                   

\begin{abstract}
In order to make the numerical simulation of atherosclerotic plaque growth feasible, a temporal homogenization approach is employed. The resulting macro-scale problem for the plaque growth can be further accelerated by using parallel time integration schemes, such as the parareal algorithm. However, the parallel scalability is dominated by the computational cost of the coarse propagator. Therefore, in this paper, an interpolation-based coarse propagator, which uses growth values from previously computed micro-scale problems, is introduced. For a simple model problem, it is shown that this approach reduces both the computational work for a single parareal iteration as well as the required number of parareal iterations.
\end{abstract}

\section{Introduction}

We are concerned with the numerical simulation of atherogenesis, that is, the growth of plaque inside arteries. For instance, in case of rupture, this can lead to heart attacks or strokes, which are often fatal. While the plaque growth takes place over months to years, the wall shear stress, which is an important driving force for plaque growth, varies within the second of a heartbeat. Therefore, a fully time-resolved fluid-structure interaction (FSI) simulation of plaque growth can easily require $\mathcal{O}(10^9)$ time steps; this is computationally infeasible.

In order to make numerical simulations of arterial plaque growth feasible, a temporal homogenization approach~\cite{FreiRichter2020} can be employed, decoupling the FSI and the plaque growth. In the resulting algorithm, a periodic FSI micro scale problem with fixed growth values has to be solved in each macro scale step of the plaque growth problem. The temporal micro-scale problem cannot be neglected. However, only few heartbeats have to be simulated to reach a periodic state; cf.~\cite{FreiRichterWick2016,Frei:2021:THN}. As a result, the computational cost is reduced significantly.

The time to solution can be further reduced by parallelization. Using parallelization in space typically results in a good parallel efficiency, however, even after reaching a strong scaling limit, the time to solution might remain high due to the large number of time steps. In order to parallelize further, a parallel time integration approach for the plaque growth problem has been introduced  in the preprint~\cite{FreiHeinlein2022}, based on the parareal algorithm~\cite{Lions:2001:PIT}. As usual in two level methods, the parallel scalability is limited by the cost of the coarse level, here, the coarse propagator. In particular, when the classical parareal approach is applied to the plaque growth problem, a periodic micro scale FSI problem has to be solved in each time step of the coarse propagator. To speed up the coarse propagator, an approach which reuses growth values from previous micro-scale computations has been introduced in~\cite{FreiHeinlein2022}, leading to an increase in the number of required parareal iterations, but a reduction in the overall time to solution.

In this paper, we introduce an interpolation-based coarse propagator, which reduces the computational cost while being much more accurate compared to the re-use approach. The new interpolation approach uses growth values from all previous micro-scale computations for the interpolation and, thus, becomes more accurate with each parareal iteration. We test it on a simple two-dimensional model, for which it even speeds up the convergence of the parareal iteration, compared to the classical parareal algorithm.

\section{Equations} \label{sec:equations}

We consider a time-dependent fluid-structure interaction system, where the fluid is modeled by the Navier--Stokes equations and the solid by the Saint Venant--Kirchhoff model. To account for the solid growth, we { incorporate} a multiplicative growth term { into} the deformation gradient{. The growth term} is motivated by typical plaque growth models{, such as those in}~\cite{YangJaegerNeussRaduRichter2015,FreiRichterWick2016,
SilvaJaegerNeussRaduSequeira2020}.

\subsection{Fluid-structure interaction}

We consider a partition of an overall domain $\Omega(t) = \FL(t) \cup \IN(t)\cup \SO(t)$ into a fluid part $\FL(t)$, an interface $\IN(t)$ and a solid part $\SO(t)$. The blood flow and its interaction with the surrounding vessel wall is modelled by the following {non-stationary} FSI system:
{\begin{equation}\label{fullplaquemodel}
	\begin{aligned}
	\rho_f(\partial_t \vt_f+\vt_f\cdot\nabla\vt_f) -
	\operatorname{div}\,\sigmat_f&= 0, &\qquad
	\operatorname{div}\,\vt_f&=0 \qquad\text{ in }\FL(t),\\
	\rho_s \partial_t\hat{\vt}_s -
	\operatorname{div}\,(\hat{\Ft}_e\hat{\Sigmat}_e) &=0, &\qquad
	\partial_t \hat\ut_s - \hat\vt_s &=0 \qquad \text{ in }\hat{S},\\
	\sigmat_f\nt_f+\sigmat_s\nt_s &=0, &\qquad
	\vt_f&=\vt_s \quad\;\; \text{ on }\IN(t).
	\end{aligned}
	\end{equation}}
Here, $\vt_f$ and $\hat{\vt}_s$ denote the fluid and solid velocit{ies}, respectively, and 
$\hat{\ut}_s$ the solid displacement.
Quantities with a ``hat'' are defined in Lagrangian coordinates, quantities without a ``hat'' in the current Eulerian coordinate framework. Two quantities $\hat{f}(\hat{x})$ and $f(x)$ correspond to each other by a $C^{1,1}$-diffeomorphism $\hat{\xi}: \hat{\Omega}\to \Omega(t)$ and the relation $\hat{f} = f\circ \hat{\xi}$. Below{,} we will also {make use of the} solid deformation gradient $\hat{F}_s = I + \hat{\nabla } \hat{u}_s$, which is the derivative of $\hat{\xi}$ in the solid part.
The constants $\rho_f$ and $\rho_s$ are the densities of blood and vessel
wall, and $\nt_f$ and $\nt_s$ are outward pointing normal vectors of the fluid and solid domain{s}, respectively.

{Furthermore, let} $\sigmat_f=\rho_f\nu_f(\nabla\vt_f+\nabla\vt_f^T)-p_f I$ denote the Cauchy stress tensor of the fluid, where $\nu_f$ is the kinematic viscosity of blood. {On the other hand,} $\sigmat_s$ denotes the Cauchy stress tensor of the solid{, which is related} to the Piola-Kirchhoff stress $\hat{\Sigmat}_e$ (which is defined below) by the Piola transformation $\hat{\sigmat}_s = \hat{J}_e^{-1} \hat{\Ft}_e\hat{\Sigmat}_e \hat{\Ft}_e^T$.

\begin{figure}
	\centering
	\resizebox*{0.65\textwidth}{!}{
	\begin{tikzpicture}
	\draw[fill=gray!30,line width=1pt] (-5,-2) -- (5,-2) -- (5,-1) -- (-5,-1) -- cycle;
	\node () at (2.5,-1.5) {$\mathcal{S}$};
	\draw[line width=1pt] (-5,-1) -- (5,-1) -- (5,1) -- (-5,1) -- cycle;
	\draw[fill=gray!30,line width=1pt] (-5,1) -- (5,1) -- (5,2) -- (-5,2) -- cycle;
	\node () at (2.5,1.5) {$\mathcal{S}$};
	
	\node () at (3.0,0.5) {$\mathcal{F}$};
	
	\draw[|<->|,line width=0.5pt] (-3.5,-2) -- (-3.5,-1); \node () at (-3.75,-1.5) {1};
	\draw[|<->|,line width=0.5pt] (-3.5,-1) -- (-3.5,0); \node () at (-3.75,-0.5) {1};
	\draw[|<->|,line width=0.5pt] (-3.5,0) -- (-3.5,1); \node () at (-3.75,0.5) {1};
	\draw[|<->|,line width=0.5pt] (-3.5,1) -- (-3.5,2); \node () at (-3.75,1.5) {1};
	
	\draw[|<->|,line width=0.5pt] (-5,-2.25) -- (0,-2.25); \node () at (-2.5,-2.5) {5};
	\draw[|<->|,line width=0.5pt] (0,-2.25) -- (5,-2.25); \node () at (2.5,-2.5) {5};
	
	\node () at (-5.5,1.5) {$\Gamma_s$};
	\node () at (5.5,1.5) {$\Gamma_s$};
	\node () at (-5.5,-1.5) {$\Gamma_s$};
	\node () at (5.5,-1.5) {$\Gamma_s$};
	
	\node () at (-5.5,0.25) {$\Gamma_f^{\text{in}}$};
	\node () at (5.5,0.25) {$\Gamma_f^{\text{out}}$};
	
	\node () at (-2.0,0.75) {$\Gamma$};
	\node () at (-2.0,-0.75) {$\Gamma$};
	
	\draw[line width=1.0pt,dotted] (-5.5,0) -- (5.5,0);
	\draw[line width=1.0pt,dotted] (0,-2.5) -- (0,2.5);
	
	\draw[fill=black] (0,-1) circle (0.1); \node () at (-0.5,0.75) {(0,1)};
	\draw[fill=black] (0,1) circle (0.1); \node () at (-0.5,-0.75) {(0,-1)};
	
	\draw[fill=black] (0,0) circle (0.05); \node () at (0.5,0.25) {(0,0)};
	\end{tikzpicture} }
	\caption{Sketch of the computational domain centered at the origin $(0,0)$; $\FL$ and $\SO$ are the fluid and solid parts, respectively, and the solid lines correspond to the fluid-solid interface $\Gamma$. The plaque growth is initiated at $(0,\pm1)$.
		\label{fig:conf}}
\end{figure}
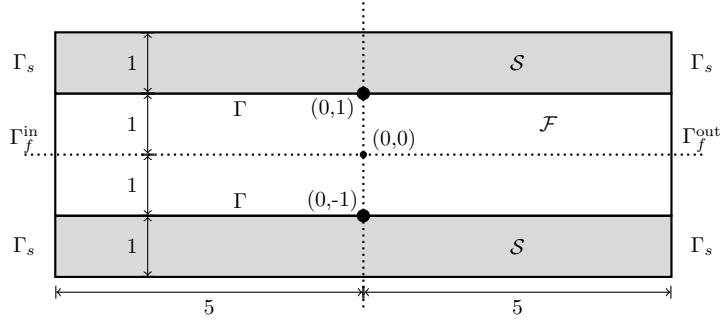

A sketch of the computational domain, which will be used in the numerical examples below, is given in~\cref{fig:conf}. We split the outer boundary of $\Omega$ into a solid part $\Gamma_s$ with homogeneous Dirichlet conditions {for the solid displacement}, a fluid part $\Gamma^\text{in}_f$ with an inflow Dirichlet condition {for the fluid velocity} and an outflow part $\Gamma^\text{out}_f$, where a \textit{do-nothing {boundary} condition} is imposed. {In particular, t}he boundary data is given by 
\begin{align}\label{problem:long-scale-boundary}
\vt_f = \vt^\text{in} \,\text{ on }\Gamma^\text{in}_f,
\quad
\rho_f\nu_f(\nt_f\cdot\nabla)\vt_f-p_f\nt_f =0
& \,\text{ on }\Gamma^\text{out}_f,\quad
\hat{\ut}_s=0\,\text{ on }\hat{\Gamma}_s,
\end{align}
where $\vt^\text{in}$ is the inflow velocity on $\Gamma^\text{in}_f$.

\subsection{A simple model for plaque growth} \label{sec:solid_growth}

Developing a realistic model {for the} plaque growth at the vessel wall is a complex task that involves the interaction of many different molecules 
and species{;} see for example~\cite{SilvaJaegerNeussRaduSequeira2020}. Furthermore, the plaque growth will also strongly depend on the geometry and material model of the arterial wall, {as well as the interaction with the blood flow}.
In this contribution, {we restrict ourselves to a simplified plaque growth model and instead focus on algorithmic aspects. 
}
{More specifically, w}e consider a simple ODE-based model  taken from~\cite{FreiRichterWick2016,Frei:2021:THN} that focuses on the influence of the concentration of foam cells $c_s$ on the growth{; the spatial distribution of the growth is prescribed in this model.} 


{Therefore, let} the spatial distribution of the plaque growth {be fixed up to a constant factor, which corresponds to the foam cell concentration $c_s$. Moreover, we assume that} the evolution of the foam cell concentration $c_s$ can be described by a simple ODE{, which depends} on the wall shear stress $\sigmat_f^{WS}$ at the vessel wall:
\begin{equation}\label{reaction}
\begin{aligned}
\partial_t c_s & = \gamma(\sigmat_f^{WS},c_s) := \frac{\alpha}{
\big(1+c_s\big)}\Bigg(1+\frac{\big\|\sigmat_f^{WS}\big\|_{L^2(\Gamma)}^2}{\sigma_0^{2}}\Bigg)^{-1}, \qquad
\sigmat_f^{WS} & := 
\big(I_d-\nt_f\nt_f^T\big) \sigmat_f
\nt_f.
\end{aligned}
\end{equation}
The reference wall shear stress $\sigma_0$ and the scale separation parameter $\alpha$ are parameters of the growth model. For cardiovascular plaque growth, we  typically have $\alpha = \mathcal{O}(10^{-7})\,\text{s}^{-1}$; {see the discussion in}~\cite{FreiRichter2020}. 

{As mentioned earlier, w}e model the solid growth by a multiplicative splitting of the 
deformation gradient $\hat \Ft_s${. In particular, let $\hat\Ft_e$ be the elastic part and $\hat\Ft_g$ a growth function, such that}
\begin{equation}\label{elasticgrowth-reference}
\hat \Ft_s=\hat\Ft_e\hat\Ft_g\quad\Leftrightarrow\quad
\hat\Ft_e = \hat\Ft_s\hat\Ft_g^{-1} =
[I+\hat\nabla\hat\ut_s]\hat\Ft_g^{-1};
\end{equation}
cf.~\cite{RodriguezHogerMcCulloch1994, YangJaegerNeussRaduRichter2015, FreiRichterWick2016}.
In the ODE model, we use the following growth function depending on $c_s$
\begin{equation}\label{num:1:growth}
\hat g(\hat x,\hat y,t) =  1+ c_s \exp\left(-\hat
x^2\right)(2-|\hat y|),\quad 
\hat\Ft_g(\hat x,\hat y,t):=\hat g(\hat x,\hat y,t)\,I.
\end{equation}
{Hence, the growth rate is influenced by the fluid-structure interaction only in form of the}
variable $c_s$. {From~\cref{num:1:growth} and the fact that the computational
domain is centered around the origin, we obtain that growth is concentrated at $(0,\pm1)$; 
see~\cref{fig:conf}.}
It follows that
\begin{equation}
\hat \Ft_g = \hat gI\quad\Rightarrow\quad
\hat\Ft_e:=\hat g^{-1} \hat\Ft_s,\label{FgFe}
\end{equation}
and the elastic Green--Lagrange strain is given by
\begin{equation}\label{Es}
\hat{\bs{E}}_e = \frac{1}{2}(\hat
{\bs{F}}_e^T\hat{\bs{F}}_e-I) = \frac{1}{2}(\hat g^{-2}\hat
{\bs{F}}_s^T\hat{\bs{F}}_s-I)
\end{equation}
resulting in the Piola--Kirchhoff stresses ({see~\cite{FreiRichterWick2016} for details})
\begin{equation}\label{Piola}
\hat\Ft_e\hat\Sigmat_e = 2\mu_s \hat{\bs{F}}_e\hat{\bs{E}}_e +
\lambda_s  \operatorname{tr}(\hat{\bs{E}}_e)\hat{\bs{F}}_e
= 2\mu_s \hat g^{-1}\hat{\bs{F}}_s\hat{\bs{E}}_e +
\lambda_s  \hat g^{-1}\operatorname{tr}(\hat{\bs{E}}_e)\hat{\bs{F}}_s.
\end{equation}

\section{Numerical framework}
\label{sec:num}

To solve the FSI problem~\cref{fullplaquemodel} we use an Arbitrary Lagrangian Eulerian (ALE) approach; see, e.g.,~\cite{RichterBuch}. {For a more detailed description of the numerical framework, we refer to~\cite{FreiHeinlein2022}.}

\subsection{Temporal two-scale approach} \label{sec:two_scale}

Even for the simplified two-dimensional configuration considered in this work, 
a resolution of the micro-scale dynamics with a scale of milliseconds to seconds is unfeasible over the complete time interval of interest $[0,T_{\text{end}}]$, with $T_{\text{end}}$ being several months 
up to a year. For instance, when considering a relatively coarse micro-scale time step of $\delta \tau = 0.02$s, the number of time steps required to simulate a time frame of a whole year would be $365 \cdot 86\,400\cdot \frac{1s}{\delta\tau} \approx 1.58\cdot 10^9$, each step corresponding to the solution of a mechano-chemical FSI problem. 

{As a remedy,}
we apply the two-scale approach of Frei and Richter {developed} in~\cite{FreiRichter2020}{,} which can be seen as a variant of the Heterogeneous Multiscale Method (HMM)~\cite{EEngquist2003}. 
{In particular, in~\cite{FreiRichter2020},} a periodic-in-time micro-scale problem is 
solved in each time step of the macro scale, for instance each day. 
The growth function $\gamma(\sigmat_f^{WS})$ is then averaged by integrating over one period of the heart beat, its average will be denoted by $\overline{\gamma}(\sigmat_f^{WS})$. This averaged growth function is applied to advance the foam cell concentration {based on}~\cref{reaction}. 



Firstly, we divide the macro-scale time interval $[0,T_{\text{end}}]$ into $N_l$ time steps of size $\delta t$
\begin{align}\label{deft}
0 = t_0 < t_1 < ... < t_{N_l} = T, \qquad N_l = \frac{T_{\text{end}}}{\delta t}.
\end{align} 
{Since} $c_s$ varies significantly on the macro scale only, using $c_s(t_m)$ as a fixed value for the growth variable on the micro scale is a reasonable approximation in the time interval $[t_m,t_{m+1}]$. Then, one cycle of the pulsating blood flow problem (around $1$\;s) is to be resolved on the micro scale 
\begin{align}\label{deftau}
0 = \tau_0 <\tau_1 <....< \tau_{N_s} = 1s, \qquad N_s = \frac{1s}{\delta \tau}
\end{align}
{to compute the averaged growth $\overline{\gamma}(\sigmat_f^{WS})$.}
It has been shown in~\cite{FreiRichter2020} (for a simplified flow configuration) that this approach leads to a model error ${\cal O}(\epsilon)$ compared to a full resolution of the micro scale, where $\epsilon=\frac{1s}{T_{\text{end}}}$ denotes the ratio between macro and micro time scale and is in the range of ${\cal O}(10^{-7})$ for a typical cardiovascular plaque growth problem. From relations~\eqref{deft} and \eqref{deftau}, we have $\delta \tau = \epsilon \delta t$.
We note that, in the model problem formulated above, scale separation is induced by the parameter $\alpha={\cal O}(\epsilon)$.

Within the two-scale algorithm, periodic micro-scale problems need to be solved. 
If {suitable} initial conditions for $\wt^0:=(\vt^0, \ut^0)$ on the micro-scale are {available}, for example, from {the micro scale problem from the} previous {macro} time step, convergence to the periodic solution may be obtained by simulating a few cycles of the micro-scale problem, due to the dissipation of the flow problem; see~\cite{FreiRichter2020,Frei:2021:THN}. {After each cycle, we can check numerically} if the solution is sufficiently close to a periodic state. {In particular, w}e apply a stopping criterion based on the computed averaged growth value:
\begin{align*}
|\overline{\gamma}(\sigmat_f^{WS,r}) - \overline{\gamma}(\sigmat_f^{WS,r-1})| < \epsilon_p,
\end{align*}
where $r=1,2,...$ denotes the iteration index with respect to the number of cycles of the micro-scale problem. The two-scale algorithm is summarized in~\cref{twoscale}.

\begin{algorithm}
\caption{Two-scale algorithm\label{twoscale}}
  Set starting values $\wt^{0,0}=(\vt^{0,0}, \ut^{0,0})$ and time step sizes $\delta t, \delta T,\, N_s=\frac{1}{\partial t}, N_l=\frac{T_{\text{end}}}{\delta T}$.
  \\
  \For{$n=1,2,\dots, N_l$}{ 
  \begin{enumerate}
    \item[1.)] \textbf{Micro problem}: Set $r\gets 0$\\
    \While{$|\overline{\gamma}(\sigmat_f^{WS,r}) - \overline{\gamma}(\sigmat_f^{WS,r-1})| > \epsilon_p$}{
    \begin{itemize}
  \item[1.a)]
  Solve micro-scale problem~\eqref{fullplaquemodel} 
    in $I_n =(t_n,t_n+ 1s)$
    \[    
    \{\wt^{r,0},c_s^{n-1}\}\mapsto 
    \{\wt^{r,m}\}_{m=1}^{N_s}
    \]
  \item[1.b)] Compute the averaged growth function
    \[
    \overline{\gamma}(\sigmat_f^{WS,r+1}) = \frac{1}{N_s}\sum_{m=1}^{N_s}
    \gamma(\sigmat_f^{WS,m}(\vt^{r+1,m}), c_s^{n-1})
    \]
    and set $\wt^{r+1,0}= \wt^{r,N_s}, \, r\gets r+1$.
    \end{itemize}}
  \item[2.)] \textbf{Macro problem}: Update the foam cell concentration $c_s^n$ by~\cref{reaction}.
    \end{enumerate}\vspace{-0.4cm}
}
\end{algorithm}

We note that the macro problem in step 2.) is {computationally} very cheap {since} it consists of one time step of an ordinary differential equation. Before each macro step, however, a micro-problem needs to be solved in step 1.a) of~\cref{twoscale}. The solution of this micro problem is typically extremely expensive, {as} a time-dependent FSI problem needs to be solved.  Considering a relatively coarse micro-scale discretization of $\delta t = 0.02$\,s, as {the one used in~\cref{sec:results}}, 
$50$ time steps are necessary to compute a single period of the heart beat. The simulation of two or more cycles is typically necessary to obtain a near-periodic state, that fulfills the stopping criterion in 1); see~\cite{FreiRichter2020,Frei:2021:THN}. In a realistic scenario, each time step of the micro problem corresponds to the solution of a complex three-dimensional FSI problem, which already makes the solution of {a single} micro problem (i.e., $\geq 100$ FSI time steps) very costly. We will thus in the following estimate the computational cost only by means of the number of micro problems that need to be solved. Even when a PDE of convection-diffusion-reaction type is considered in the plaque growth model {instead}, the cost of the micro problems dominates the overall {computational work} by far{; see the discussion and the numerical results in~\cite{FreiHeinlein2022}.}

\section{Parallel time stepping}
\label{sec:parallel}

In the preprint~\cite{FreiHeinlein2022}, we have introduced parallel time stepping algorithms {for the macro problem} based on the parareal algorithm to further reduce the computational times. It turns out that the coarse propagator limits the maximum possible speed-up severely. A first algorithm to reduce the cost of the coarse propagator {by re-using growth values from previous micro-scale computations} has already been proposed in~\cite{FreiHeinlein2022}{; see also~\cref{sec.reusage}.} In this section, we will first review the developments in~\cite{FreiHeinlein2022} and then propose a {different approach} to reduce the cost of the coarse propagator, based on an interpolation of the previously computed growth values.



\subsection{The parareal algorithm}

First, the time interval of interest $[0,T_{\text{end}}]$ on the macro scale is divided into $P$ coarse sub-intervals $I_p=[T_{p-1}, T_p]$ of equal size. {In the setting of parallelization in time, the sub-intervals} are distributed among {$P$ parallel} processes $p=1,...,P$:
\begin{align}\label{eq:coarse_disc}
0 = T_0 < T_1 < ... < T_P = T_{\text{end}}.
\end{align}

In order to define the parareal algorithm, suitable fine- and coarse problems need to be introduced. {Following the work~\cite{FreiHeinlein2022}}, we will apply the parareal algorithm only on the macro scale; hence, both the fine and the coarse scale of {the parareal algorithm} correspond to the macro scale of the homogenization approach. The fine problem advances the growth variable $c_s$ from time $T_p$ to $T_{p+1}$ by applying the two-scale algorithm (\cref{twoscale}) with a fine time step size $\delta t$ (e.g.\,$0.3$ days) on the corresponding fine
time discretization of {the interval} $[T_p, T_{p+1}]$:
\begin{align*}
T_p = t_{p,0} < t_{p,1} < ... < t_{p,N_p} = T_{p+1}, \qquad N_p = \frac{T_{p+1}-T_p}{\delta t}, \qquad  t_{p,q} := t_{p\cdot N_p +q} = t_{p,q-1} + \delta t.
\end{align*}
We use capital letters $T_p$ to denote the coarse discretization of $[0,T_{\text{end}}]$ into $P$ parts corresponding to the processes and small letters $t_i$ to denote the finer discretization on the {parallel} processes; the two discretizations yield the first level (fine problem) and second level (coarse problem) of the parareal algorithm. 
In the example studied in~\cref{sec:results}, the time interval has length $T_{\text{end}}=300$ days, $\delta T$ is 30 days for $P=10$, while $\delta t$ is chosen as 0.3 days. On the micro scale, the times $\tau_i$ and time step size $\delta\tau$ are defined locally in $[t_i, t_i + 1s]$. {In particular,} 
we will use $\delta \tau =0.02$\,s {here}.
Note that the micro scale influences the parareal algorithm only indirectly {through} the temporal two-scale approach.


On each process $p=1,\ldots,P$, the fine propagator $p$ consists of a time stepping procedure to advance $c_s(T_p)$ to $c_s(T_{p+1})$ with the fine {time step size} $\delta t$. We abbreviate this computation by
\begin{align*}
c_s(T_{p+1}) = {\cal F}(c_s(T_p)).
\end{align*}
{In order to obtain parallel scalability, t}he coarse propagator needs to be much cheaper, since it advances the foam cell concentration globally on $[0,T_{\text{end}}]$ and thus introduces synchronization.
A natural choice 
is to use the same two-scale algorithm as for the micro problems, but with a much larger time-step size $\delta T=T_{p+1}-T_p$. 
This means that in total $P$ macro steps are required to solve the coarse problem on $[0,T_{\text{end}}]$, which includes the solution of $P$ micro-scale FSI problems. We introduce the notation
\begin{align*}
c_s(T_{p+1}) = {\cal C}(c_s(T_p))
\end{align*}
for one time step of this coarse propagator.

Then, given an iterate $\overline{c}_s^{(k)}$ for some $k\geq 0$, the parareal algorithm computes $\overline{c}_s^{k+1}$ by setting
\begin{align}\label{eq:parareal}
\overline{c}_s^{(k+1)}(T_{p+1}) = {\cal C} (\overline{c}_s^{(k+1)}(T_p)) +{\cal F}( \overline{c}_s^{(k)}(T_p)) - {\cal C}(\overline{c}_s^{(k)}(T_p)) \quad \text{ for }\, p=0,...,P-1.
\end{align}
As the fine-scale contributions depend only on the previous iterate $\overline{c}_s^{(k)}$ the contributions from each process $p=1,\ldots,P$ can be computed in parallel.


Let us discuss the computational cost of the parareal algorithm in more detail. As mentioned above, {the dominant part of the}
computational cost in the two-scale algorithm (\cref{twoscale}) {corresponds to} the solution of the micro problems.
The first term in~\cref{eq:parareal} requires the solution of one micro problem in each coarse time step $T_p\to T_{p+1}$. 
Within the fine-scale propagator (second term in~\cref{eq:parareal}) a maximum of $N_p=\lceil N_l/P\rceil$ time steps need to be computed per process, where $\lceil f \rceil$ denotes the next-biggest natural number to $f$ ($\widehat{=}$ \texttt{ceil(f)}). Each time step comes with the solution of one micro problem. The last term in \cref{eq:parareal} has already been computed in the previous iteration (compare the first term on the right-hand side of the same equation) and therefore introduces no additional computational cost. Considering that all  fine-scale problems ($p=1,\ldots,P$) can be solved in parallel and {assuming} that each micro problem requires approximately the same time, {the computing time for} $k_{\text{par}}$ iterations of the parareal algorithm {on $P$ processes corresponds to} 
\begin{align}\label{compcost}
	\underbrace{k_{\text{par}}\cdot \lceil N_l/P \rceil}_{\text{fine level ($P$ parallel processes)}} + \underbrace{(k_{\text{par}}+1) \cdot P}_{\text{coarse level (1 serial process)}}
\end{align}
micro problems in a serial computation.

\subsection{Re-Usage of Growth Values}
\label{sec.reusage}

{Once the number of processes exceeds  $P\geq \sqrt{N_l}$, the cost of the coarse-scale propagator gets dominant; cf.~\cref{compcost}.} In~\cite{FreiHeinlein2022}, we propose an alternative coarse-scale propagator which avoid{s} the solution of micro problems. Instead, we use the values $\overline{\gamma}(\sigmat_f^{WS}(t_{p,i}))$ computed within the fine-scale propagators. For this purpose, all values $\overline{\gamma}_{p\cdot N_p + i}:=\overline{\gamma}(\sigmat_f^{WS}(t_{p,i}))$ computed on the fine scale on all processes $p=1,...,P$ are stored
for all time steps $i=1,...,N_p$ and re-used within the coarse propagator.

{Since} no new micro problems need to be solved and approximations of the growth values $\overline{\gamma}_j$ are available for all fine time steps $j=1,...,N_l$, the coarse propagator can even {be applied with the fine-scale} 
time step {size} $\delta t$. The only additional cost is to advance the foam cell concentration by the ODE~\eqref{reaction}, which is clearly negligible compared to the solution of the micro problems in the fine-scale propagator. The only coarse-scale propagator, which can not be replaced, is the first one which is used for initialization. The total runtime for $k_{\text{par}}$ iterations of the re-usage variant is
\begin{align}\label{cc_reusage}
	\underbrace{
	k_{\text{par}} \cdot \lceil N_l/P \rceil}_{\text{fine level ($P$ parallel processes)}} + \underbrace{P}_{\text{coarse level (1 serial process)}}
\end{align}
serial micro problems{. Compared to~\cref{compcost}, this approach saves the computation of $k_{\text{par}} P$ micro problems on the coarse scale; however, since it uses an approximation for the growth values, the number of parareal iterations $k_{\text{par}}$ might be slightly higher.} For details on the resulting algorithm, we refer to~\cite{FreiHeinlein2022}.

%

\section{Interpolation Approach} \label{sec:interpolation}

The drawback of the re-use approach is that the growth values $\overline{\gamma}(\sigmat_f^{WS})(t_{j})$ {are} computed based on the previous iterate of the foam cell concentration $c_s^{(k)}(t_{j})$ instead of the parareal iterate $\overline{c}_s^{(k+1)}(t_{j})$, which would be the most accurate approximation available in the coarse propagator. While we can expect that these are quite close after a few iterates, we have observed in~\cite{FreiHeinlein2022} that this is not necessarily the case in the very first parareal iterate.
Thus, typically a few more parareal iterations are required to reach a stopping criterion.

{Here, we propose another approach, which may lead to better approximations of the growth values $\overline{\gamma}(\sigmat_f^{WS})(t_{j})$. In particular, we} 
store again the pairs of computed foam cell concentrations and growth values $(c_s^{(k)}(t_{i}), \overline{\gamma}(\sigmat_f^{WS}(t_{i})))$ over all iterates $k$. In contrast to the re-usage variant, we 
interpolate {linearly} from the {two} values which are closest to $\overline{c}_s^{(k+1)}(t_{j})$. For this purpose, we create an ordered list of values $\{c_s^i\}_{i=1}^N$ and their corresponding growth values $\{\gamma_{\text{int}}^i\}_{i=1}^N$ \textit{on the fly} 
{during the parareal iterations}
and compute the linear interpolation from the values that fulfill $c_s^i < \overline{c}_s^{(k+1)}(t_{j}) < c_s^{i+1}$ when this is needed:
\begin{align}\label{eq:interpol}
\overline{\gamma}_{\text{int}}(\sigmat_f^{WS})(t_{j}) := \frac{c_s^{i+1} - \overline{c}_s^{(k+1)}(t_{j})}{c_s^{i+1} - c_s^i} \gamma_{\text{int}}^{i} \,+\, \frac{\overline{c}_s^{(k+1)}(t_{j}) - c_s^i}{c_s^{i+1} - c_s^i} \gamma_{\text{int}}^{i+1};
\end{align}
The resulting algorithm is given in~\cref{alg:interpol}. As in the re-usage variant no micro-problems need to be solved within the coarse-scale propagator. The total runtime corresponds thus again to the serial solution of $k_{\text{par}}\cdot \lceil N_l/P\rceil + P$ micro problems and the optimal choice in terms of speed-up is again $P \approx \sqrt{k_{\text{par}} N_l}$; see the discussion in Section 4.2.2 in~\cite{FreiHeinlein2022}. {However, compared to the re-usage algorithm, we expect that fewer iterations of the parareal algorithm $k_{\text{par}}$ are required to fulfill the stopping criterion.}

{This approach can be regarded as a simple linear reduced order model for the micro problems. Since the interpolation data is updated during the parareal iterations, the model predictions may actually improve from iteration to iteration.}

\begin{algorithm}
	\caption{Parallel time stepping with interpolation of growth values\label{alg:interpol}}
	\begin{enumerate}
		\item[(I)] \textbf{Initialization}: Compute $\big\{(\overline{c}_s^{(0)}(T_p), \overline{\wt}^{(0)}(T_p)\big\}_{p=1}^P$ by means of~\cref{twoscale}
		with a coarse macro time step $\delta t := T_{p+1}- T_p$ on the master process. Set $k\gets 0$\\
		\item[(II)]\While{$|c_s^{(k+1)}(T_{\text{end}}) - c_s^{(k)}(T_{\text{end}})| > \epsilon_{\text{par}}$}{
			\begin{itemize}
				\item [(II.a)] \textbf{Fine problem} on each process $p=1,...,P$: 
				\begin{itemize}
					\item[(i)] Initialize $c_s^{(k+1)}(T_p)=\overline{c}_s^{(k)}(T_p), \,\wt^{(k+1)}(T_p)  = \wt^{(k)}(t_{p-1, N_{p}})$
					\item[(ii)] Compute $\big\{(c_s^{(k+1)}(t_{p,q}), \wt^{(k+1)}(t_{p,q})\big\}_{q=1}^{N_p}$ by \cref{twoscale} with fine time step $\delta t$
					\item[(iii)] Add the pair $(c_s^{(k+1)}(t_{p,q}), \overline{\gamma}(\sigmat_f^{WS})(t_{p,q}))$ to the list of growth value pairs and order it by the value of its first component
				\end{itemize}  
				\item[(II.b)] \textbf{Coarse problem} on the master process:\\
				\For{$j=0,...,N_l-1$}{
					\begin{itemize}
						\item[(i)] Compute the interpolated growth value $\overline{\gamma}_{\text{int}}(t_j)$ corresponding to $\overline{c}_s^{(k+1)}(t_{j})$ by~\cref{eq:interpol} 
						\item[(ii)] Advance
						$\overline{c}_s^{(k+1)}(t_{j+1}) = \overline{c}_s^{(k+1)}(t_{j}) + \delta t\, \overline{\gamma}_{\text{int}}(t_j)$.
					\end{itemize}
				}\vspace{-0.3cm}
		\end{itemize}}\vspace{-0.55cm}
	\end{enumerate}
\end{algorithm}

\section{Numerical Results} \label{sec:results}


We show numerical results for a numerical example taken from~\cite{FreiHeinlein2022};
see~\cref{fig:conf} {for a visualization of the configuration}. The material parameters are chosen close to {those of} a realistic plaque growth problem; for details {on the choice of the parameters,} we refer to~\cite{FreiHeinlein2022}.
As an inflow boundary condition, a pulsating
velocity inflow profile is prescribed on $\Gamma_f^\text{in}$ .

\begin{figure}
\begin{center}
\includegraphics[width=0.7\textwidth]{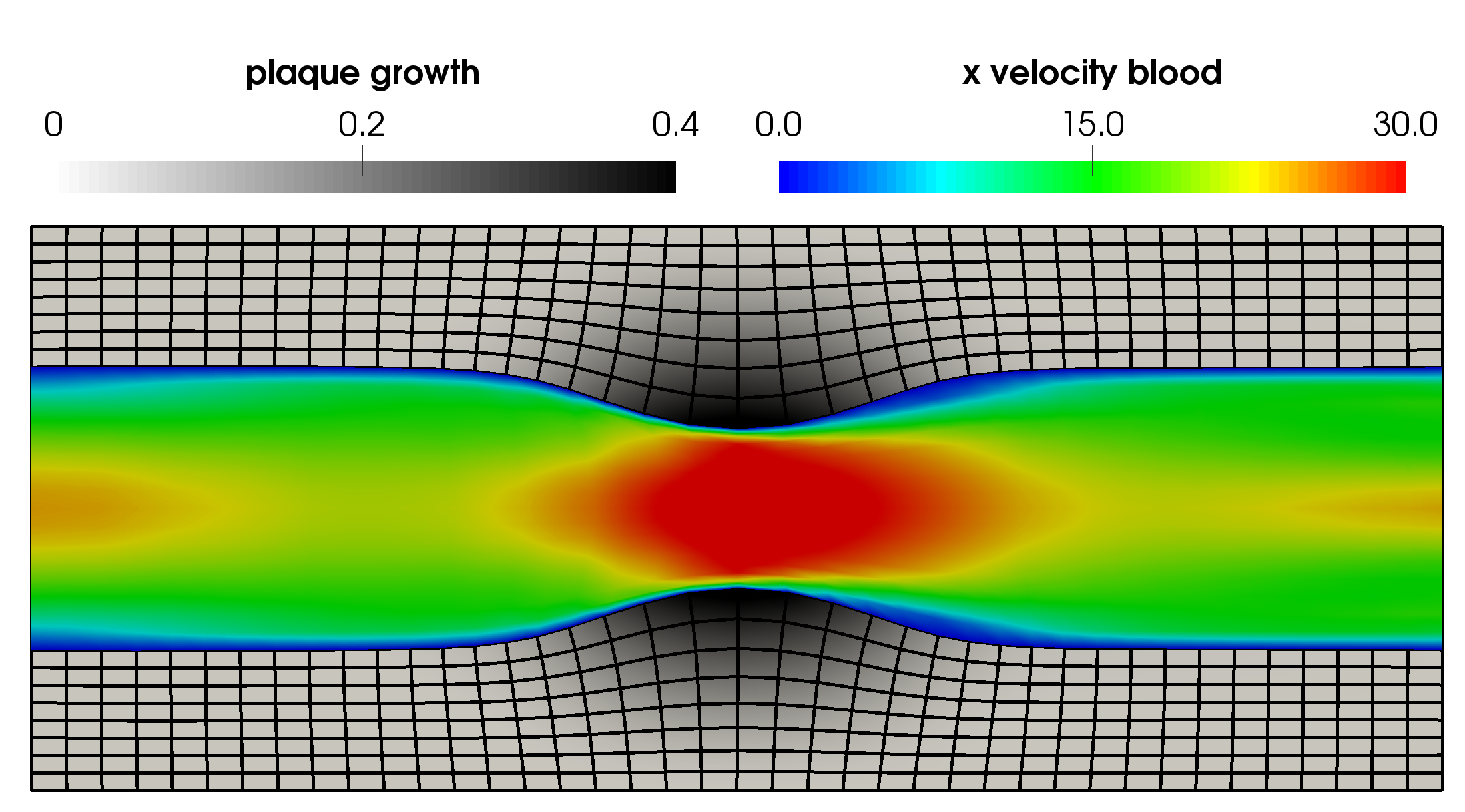}
\end{center}
\caption{\label{fig:visu} Visualization of the plaque growth at time $t=300$ days. The horizontal velocity (in cm/s) and the vertical displacement (in cm) are shown on the deformed domain at micro time $\tau=0.34$\,s.}
\end{figure}

For time discretization, we use a backward Euler method for the FSI problem
and a forward Euler scheme for the ODE growth model.
For spatial discretization, we use biquadratic ($Q_2$) equal-order finite elements for all variables and LPS stabilization~\cite{BeckerBraack2001} for the fluid problem. 
Our mesh, containing both fluid and solid, consists of 160 rectangular grid cells; this corresponds to {relatively small FSI problem with} a total of 3\,157 degrees of freedom. 
The time step sizes are chosen as $\delta \tau=0.02$\,s and $\delta t=0.3$ days (i.e., $N_l=1\,000$); the tolerance for periodicity of the micro-scale problem {is chosen} as $\epsilon_p=10^{-3}$. All the computational results have been obtained with the finite element library Gascoigne3d~\cite{Gascoigne} using a fully monolithic approach for the FSI problem following Frei, Richter \& Wick~\cite{FreiRichterWick2016, RichterBuch}. {All computations are performed in serial, and the parallel performance is discussed based on the computing time measured in terms of serial micro problems, as discussed in~\cref{sec:parallel}.}
A visualization of the plaque growth at the end time $T_{\text{end}}$ is shown in Figure~\ref{fig:visu}.

In Tables~\ref{tab:csk},\,\ref{tab:para2} and~\ref{tab:interpol}, we show results for the standard parareal {algorithm as well as the re-usage and the interpolation variants}, respectively. In all cases, the stopping criterion 
\begin{align} \label{eq:stop}
|c_s^{(k+1)}(T_{\text{end}}) - c_s^{(k)}(T_{\text{end}})| < \epsilon_{\text{par}}=10^{-3}
\end{align}
is used and numerical results {for varying numbers of parallel processes $P$} 
are provided.

\begin{table}[t]
\begin{center} \small
\begin{tabular}{|l|rrrrr|r|}
\hline
k &$P=10$ &$P=20$ &$P=30$ &$P=40$ &$P=50$ & ref. (serial) \\
\hline
1 &$2.21\cdot 10^{-2}$ &$1.19\cdot 10^{-2}$ &$8.12\cdot 10^{-3}$
&$5.83\cdot 10^{-3}$ &$5.46\cdot 10^{-3}$ &0.63831273\\
2 &$2.24\cdot 10^{-3}$ &$5.63\cdot 10^{-4}$ &$2.62\cdot 10^{-4}$ &$1.37\cdot 10^{-4}$ &$8.61\cdot 10^{-5}$ &-\\
3 &$1.42\cdot 10^{-4}$ &$2.02\cdot 10^{-5}$&$6.53\cdot 10^{-6}$ &$2.32\cdot 10^{-6}$ &$8.30\cdot 10^{-7}$ &-\\
4 &$5.76\cdot 10^{-6}$ &- &- &- &- &-\\
\hline 
\hline 
\# mp & 450 & 230 & {\bf 222} & 235 & 260 & 1\,000 \\
speedup & 2.2 & 4.3 & {\bf 4.5} & 4.3 & 3.8 & 1.0 \\
efficiency & {\bf 22\,\%} & {\bf 22\,\%} & 15\,\% & 11\,\% & 8\,\% & 100\,\% \\
\hline
\end{tabular}
\end{center}
\caption{\label{tab:csk} Errors $|c_s^{(k)}(T_{\text{end}})-c_s^{*}(T_{\text{end}})|$ for $P=10,\ldots, 50$ for the standard parareal algorithm and a reference computation $c_s^{*}$. The time measure in terms of the number of serial micro problems (\# mp) as well as speedup and efficiency compared to the reference computation (right column) are shown; best numbers marked in \textbf{bold face}.}
\end{table}

\begin{table}[t]
\begin{center}
\centering
\setlength{\tabcolsep}{4pt}\renewcommand{\arraystretch}{1.1} \small
\begin{tabular}{|l|rrrrrrr|r|}
\hline
k &$P=10$ &$P=20$ &$P=30$ &$P=40$ &$P=50$ &$P=60$ &$P=70$ & ref. \\
\hline                                            
1 &\small $2.56\cdot 10^{-2}$ &\small $1.10\cdot 10^{-2}$ &\small $7.48\cdot 10^{-3}$ &\small $5.46\cdot 10^{-3}$ &\small $4.23\cdot 10^{-3}$ &\small$3.63\cdot 10^{-3}$ &\small $3.12\cdot 10^{-3}$ &\small 0.63831273\\
2 &\small $7.78\cdot 10^{-3}$ &\small $4.48\cdot 10^{-3}$ &\small $3.16\cdot 10^{-3}$ &\small $2.36\cdot 10^{-3}$ &\small $1.99\cdot 10^{-3}$ &\small $1.61\cdot 10^{-3}$ &\small $1.41\cdot 10^{-3}$ &\small-\\
3&\small $1.73\cdot 10^{-3}$ &\small $1.20\cdot 10^{-3}$ &\small $9.02\cdot 10^{-4}$ 
&\small $6.92\cdot 10^{-4}$ &\small $6.34\cdot 10^{-4}$ &\small $5.14\cdot 10^{-4}$ &\small $4.32\cdot 10^{-4}$ &\small-\\
4 &\small $2.32\cdot 10^{-4}$ &\small $2.28\cdot 10^{-4}$ &\small $1.87\cdot 10^{-4}$ &\small $1.31\cdot 10^{-4}$ &\small $1.26\cdot 10^{-4}$ &\small $1.10\cdot 10^{-4}$ &\small $9.98\cdot 10^{-5}$  &\small-\\
5 &\small $2.19\cdot 10^{-5}$ &\small $4.71\cdot 10^{-5}$ &\small $3.05\cdot 10^{-5}$ &\small - &\small - &\small - &\small - &\small -\\
\hline
\hline 
\# mp 	&510 &270 &200 &140 &130 &{\bf 128} &130 &1\,000 \\
speedup &2.0 &3.7 &5.0 &7.1 &7.7 & {\bf 7.8} &7.7 &1.0 \\
eff. 	&{\bf 20\,\%} &19\,\% &17\,\% &18\,\% &15\,\% &13\,\% &11\,\% &100\,\%  \\
\hline
\end{tabular}
\end{center}
\caption{\label{tab:para2} 
Errors $|c_s^{(k)}(T_{\text{end}})-c_s^{*}(T_{\text{end}})|$ for $P=10, \ldots, 70$ for the re-usage variant, best numbers marked in \textbf{bold face}.
}
\end{table}

\begin{table}[t]
	\begin{center} \small
		\begin{tabular}{|l|rrrrrr|r|}
			\hline
			k &$P=10$ &$P=20$ &$P=30$ &$P=40$ &$P=50$ &$P=60$ & ref. ($P=1\,000$) \\
			\hline
			1 &$8.3\cdot 10^{-9}$ &$3.0\cdot 10^{-8}$ &$2.5\cdot 10^{-8}$ &$9.9\cdot 10^{-9}$ &$3.6\cdot 10^{-9}$ &$5.2\cdot 10^{-9}$ &0.6383127317\\ 
			2 &$3.0\cdot 10^{-10}$ &$4.3\cdot 10^{-9}$ &$1.6\cdot 10^{-9}$ &$8.0\cdot 10^{-10}$ &$1.2\cdot 10^{-9}$ &$1.7\cdot 10^{-9}$ &-\\
			\hline \hline
			\# mp &210 &120 &98 &\textbf{90} &\textbf{90} &94 &1\,000 \\
			speedup &4.8 &8.3 &10.2 &\textbf{11.1} &\textbf{11.1} &10.6 &1.0 \\
			efficiency &{\bf 48\,\%} &42\,\% &34\,\% &28\,\% &22\,\% &18\,\% &100\,\% \\
			\hline
		\end{tabular}
	\end{center}
	\caption{\label{tab:interpol} Errors $|\overline{c}_s^{(k)}(T_{\text{end}})-c_s^{*}(T_{\text{end}})|$ for $P=10, \ldots, 60$ for the interpolation approach; best numbers marked in \textbf{bold face}.}
\end{table}

We observe that the standard parareal algorithm (\cref{tab:csk}) needs 3-4 parareal iterations to fulfill the stopping criterion. The re-usage approach (\cref{tab:para2}) needs 1-2 additional iterations, {as a result of the} 
approximation 
used in the coarse-scale propagator. Interestingly, for the interpolation variant, the values $c_s^{(1)}(T_{\text{end}})$ {on all processes agree with the reference computation in the first 8 digits already after only a single parareal iteration (the reference computation is performed with a serial time stepping).}
A second iterate is only computed, as the stopping criterion is based on a comparison between two subsequent iterates. 

This might be surprising at first sight, as the interpolation in the coarse propagator is in principle also an approximation. On the other hand, the time step size $\delta t\ll \delta T$ is much smaller compared to the standard parareal coarse propagator. In fact, the interpolation between the values from the coarse initialisation (step (I) of algorithm~\ref{alg:interpol}) yields very accurate growth values in this example, such that the coarse propagator is both much cheaper and much more accurate compared to the coarse problem in the standard parareal algorithm. {Of course, this might be partly attributed to the simple nature of the FSI and plaque growth problems. However, it shows that a reduced order model can be a powerful alternative to the standard parareal coarse propagator.}

Due to the lower number of parareal iterations, the number of micro problems to be solved in the interpolation approach is significantly smaller compared to the other approaches {for all processes}.
The minimal number of micro problems is attained for $P=40$ and $P=50$, where {only} 90 micro problems ($2\cdot25+40$ resp. $2\cdot 20+50$) need to be solved. {This} corresponds to {an expected} speed up of 11.1 compared to a serial computation. As discussed above{,} the optimal number of processes is $P \approx \sqrt{k_{\text{par}} N_l} \approx 44.72$.
If one would stop the computation already after the first iterate, which is already accurate up to 8 digits, the computational cost would reduce further to 65 (for $P=40$) resp.\,70 (for $P=50$) micro problems{; however, this would require the derivation of an error estimator instead of the stopping criterion~\cref{eq:stop}.}

The efficiency, defined as ratio between speedup and number of processes, is also significantly improved compared to the previous approaches. For $P=10$, we obtain a parallel efficiency of $48\,\%$, and for $P=40$, where a speedup of $\approx 11$ is obtained, the efficiency is still $28\,\%$. For the other two approaches the corresponding efficiencies are $20\%$ resp. $22\%$ and $11\%$ resp.\,$18\%$.

\section{Conclusion}

We have introduced an efficient coarse propagator for a scalar plaque growth model based on an interpolation of previously computed growth values. {For} 
the simple model example considered in this work, the modified parareal algorithm converges extremely fast against a reference value. Obviously, the performance needs to be confirmed in more complex examples in future.

A realistic plaque growth model typically consists of a system of partial-differential equations. In this case, an interpolation of the map $c_s \to \overline{\gamma}(c_s)$ is not straight-forward, as both $c_s$ and $\overline{\gamma}(c_s)$ are spatially distributed quantities. In this case, more sophisticated 
reduced order models or neural networks {may be used} to approximate this map. {This will be the subject of future research.}

\bibliographystyle{plain}


\begin{thebibliography}{10}

\bibitem{BeckerBraack2001}
R.~Becker and M.~Braack.
\newblock A finite element pressure gradient stabilization for the {S}tokes
  equations based on local projections.
\newblock {\em Calcolo}, 38(4):173--199, 2001.

\bibitem{Gascoigne}
R.~Becker, M.~Braack, D.~Meidner, T.~Richter, and B.~Vexler.
\newblock The finite element toolkit {G}ascoigne3d.
\newblock http://www.gascoigne.de.

\bibitem{EEngquist2003}
W.~E and B.~Engquist.
\newblock The heterogenous multiscale method.
\newblock {\em Comm. Math. Sci.}, 1(1):87--132, 2003.

\bibitem{FreiHeinlein2022}
S~Frei and A~Heinlein.
\newblock Towards parallel time-stepping for the numerical simulation of
  atherosclerotic plaque growth.
\newblock {\em arXiv e-print: https://arxiv.org/abs/2203.06526}, 2022.

\bibitem{Frei:2021:THN}
S~Frei, A~Heinlein, and T~Richter.
\newblock On temporal homogenization in the numerical simulation of
  atherosclerotic plaque growth.
\newblock {\em PAMM}, 21(1):e202100055, 2021.

\bibitem{FreiRichter2020}
S.~Frei and T.~Richter.
\newblock Efficient approximation of flow problems with multiple scales in
  time.
\newblock {\em SIAM J Multiscale Model Simul}, 18(2):942--969, 2020.

\bibitem{FreiRichterWick2016}
S.~Frei, T.~Richter, and T.~Wick.
\newblock Long-term simulation of large deformation, mechano-chemical
  fluid-structure interactions in {ALE} and fully {E}ulerian coordinates.
\newblock {\em J Comput Phys}, 321:874 -- 891, 2016.

\bibitem{Lions:2001:PIT}
J.-L. Lions, Y.~Maday, and G.~Turinici.
\newblock {A "parareal" in time discretization of {PDE}'s}.
\newblock {\em Comptes Rendus de l'Acad\'emie des Sciences - Series I -
  Mathematics}, 332:661--668, 2001.

\bibitem{RichterBuch}
T.~Richter.
\newblock {\em Finite Elements for Fluid-Structure Interactions. Models,
  Analysis and Finite Elements.}, volume 118 of {\em Lecture Notes in Comput
  Sci and Eng}.
\newblock Springer, 2017.

\bibitem{RodriguezHogerMcCulloch1994}
E.K. Rodriguez, A.~Hoger, and A.~D. McCulloch.
\newblock Stress-dependent finite growth in soft elastic tissues.
\newblock {\em J. Biomech.}, 4:455--467, 1994.

\bibitem{SilvaJaegerNeussRaduSequeira2020}
T.~Silva, W.~J\"ager, M.~Neuss-Radu, and A.~Sequeira.
\newblock Modeling of the early stage of atherosclerosis with emphasis on the
  regulation of the endothelial permeability.
\newblock {\em J Theor Biol}, 496:110229, 2020.

\bibitem{YangJaegerNeussRaduRichter2015}
Y.~Yang, W.~J{\"a}ger, M.~Neuss-Radu, and T.~Richter.
\newblock Mathematical modeling and simulation of the evolution of plaques in
  blood vessels.
\newblock {\em J Math Biol}, pages 1--24, 2014.

\end{thebibliography}

\end{document}